\documentstyle[12pt]{amsart}

\makeatletter                                            
\renewcommand{\@begintheorem}[2]{                        
\rm \trivlist \item [\hskip \labelsep {\bf #2\ \ #1.}]   
                                }                        
\makeatother                                             

\newcommand{\newsubsection}%
{{\noindent\bf\refstepcounter{subsection}\thesubsection\ \ }}

\newcommand{\newsubsubsection}%
{{\bf\refstepcounter{subsubsection}\thesubsubsection\ \ }}

\newcommand{\ts}{\vspace{\baselineskip}\noindent{\bf Proof.}$\;\;$}
\newcommand{\ZZ}{{\bf Z}}
\newcommand{\QQ}{{\bf Q}}
\newcommand{\RR}{{\bf R}}

\newcommand{\CC}{{\bf C}}

\newcommand{\FF}{{\bf F}}

\newcommand{\PP}{{\bf P}}

\newcommand{\cC}{{\cal C}}

\newcommand{\cM}{{\cal M}}

\newcommand{\cO}{{\cal O}}

\newcommand{\tT}{{\tilde T}}

\def\DynkinEEE#1#2#3#4#5#6#7
   {\begin{picture}(124, 50)%
   \put(5, 35){\circle{4}}%
   \put(7, 35){\line(1, 0){28}}%
   \put(37, 35){\circle{4}}%
   \put(39, 35){\line(1, 0){28}}%
   \put(69, 35){\circle{4}}%
   \put(71, 35){\line(1, 0){28}}%
   \put(101, 35){\circle{4}}%
   \put(103, 35){\line(1, 0){28}}%
   \put(133, 35){\circle{4}}%
   \put(69, 33){\line(0, -11){28}}%
   \put(69, 3){\circle{4}}%
   \put(5, 41){\makebox(0, 0)[b]{$#1$}}%
   \put(37, 41){\makebox(0, 0)[b]{$#3$}}%
   \put(69, 41){\makebox(0, 0)[b]{$#4$}}%
   \put(101, 41){\makebox(0, 0)[b]{$#5$}}%
   \put(133, 41){\makebox(0, 0)[b]{$#6$}}%
   \put(69, 0){\makebox(0, 0)[t]{$#2$}}%
   \put(-60, 35){\makebox(0, 0)[lb]{$#7$}}%
\end{picture}}
\begin{document}

\title{A linear system on Naruki's  moduli space of marked cubic surfaces}
\author{Bert van Geemen}
\address{Dipartimento di Matematica, Universit\`a di Pavia,
via Ferrata 1, I-27100 Pavia, Italy}
\email{geemen@@dragon.ian.pv.cnr.it}

\begin{abstract}
Allcock and Freitag recently showed that the moduli space of marked cubic surfaces is a subvariety of a nine dimensional projective space which is defined by cubic equations. 
They used the theory of automorphic forms on ball quotients to obtain these results. Here we describe the same embedding using Naruki's toric model of the moduli space. We also give an explicit parametrization of the tritangent divisors, we discuss another way to find equations for the image and we show that the moduli space maps, with degree at least ten, onto the unique quintic hypersurface in a five dimensional projective space which is invariant under the action of the Weyl group of the root system $E_6$.
\end{abstract}

\maketitle

{\bf Introduction}

\

Recently Allcock, Carlson and Toledo \cite{ACT2} studied the moduli space of 
smooth cubic surfaces using the intermediate jacobian of the cubic threefold 
which is the triple cover of projective three space branched along a cubic 
surface. They show that this moduli space, as well as 
the moduli space of marked cubic surfaces $\cM^0$ (that is, cubic surfaces with 
an
ordered set of six skew lines) are open subsets of certain 4-ball quotients.  
The Weyl group $W(E_6)$ of the root system $E_6$ acts on $\cM^0$ 
by permuting the markings on any given cubic surface, the quotient variety is 
the moduli space of cubic surfaces. The quasi projective variety $\cM^0$ has a 
natural compactification $\cM$ given by geometrical invariant theory.
The projective variety $\cM$ coincides with the Baily-Borel compactification 
of the ball quotient.
The action of $W(E_6)$ extends to $\cM$.

Using Borcherds' work on automorphic forms on ball quotients, 
Allcock and Freitag \cite{Fr} found a
$W(E_6)$-equivariant embedding of $\cM$ in a nine dimensional projective space.
The action of $W(E_6)$ on the projective space is obtained from the unique ten 
dimensional irreducible linear representation of $W(E_6)$.
This map actually already appears in a paper by 
A.\ B.\ Coble published in 1917 \cite{Co3} (and see also \cite{Y}) where 
$\cM$ is identified with the moduli space of six points in the projective plane.  The same embedding of $\cM$ was also found by 
Matsumoto and 
Terasoma \cite{MT} who used the theta constants associated to the intermediate 
jacobians.

An explicit smooth projective compactification $\cC$ (`the cross ratio 
variety') of the moduli space $\cM$ with a 
biregular action of the Weyl group was constructed by Naruki \cite{Naruki}.
It is a modification of a toric variety associated to the root system $D_4$.
Naruki constructs and studies his model as a subvariety of the product of $270$  projective lines, each component of this map is given by a cross ratio
(of certain tritangent planes containing a given line on the cubic surface).
The Weyl group acts via permutations of these $270$ projective lines.

In this paper we explicitly identify the nine dimensional linear system on 
Naruki's model $\cC$ which defines the map $F$ to $\PP^9$ discovered by
Coble, Allcock and Freitag  (see Theorem \ref{main})
$$
F:\cC\longrightarrow \cM\quad (\subset\PP^9).
$$
We also give explicit formulas for the $W(E_6)$-action on this linear system
in section \ref{we6V}.

A tritangent plane of a cubic surface is a plane which cuts out three lines on 
the surface. If these three lines meet in a point, that point is called an 
Eckart point.   
We obtain a nice parametrization, equivariant for the Weyl group of the root 
system $F_4$, of the 45 divisors in $\cM$ which parametrize 
marked cubic surfaces with an Eckart point, see Theorem \ref{tritpar}.
A study of the linear relations between tritangent planes leads to the 
discovery that $\cM$ is the singular locus of a variety $X$ defined by six 
quintic polynomials, see \ref{quints}. The group $W(E_6)$ acts on $X$ and it would be very interesting to have a moduli interpretation for $X$.
 
The Weyl group of $E_6$ is defined as a reflection group on a real six 
dimensional vector space. Complexifying and projectivizing this vector space 
one obtains a biregular action of $W(E_6)$ on a $\PP^5$.
In his book \cite{H}, Bruce Hunt suggested an identification of the moduli
space with the unique
$W(E_6)$-invariant quintic hypersurface $I_5$ in $\PP^5$. In section 
\ref{6dire}
we construct a dominant rational map $\Sigma:\cM\longrightarrow \PP^5$ 
which is equivariant for the action of $W(E_6)$ and we show that its image is 
$I_5$ (Thm.\ \ref{i5}), but, unfortunately, this map has degree at least 10
(Thm \ref{thmdeg}).

The results of this paper are obtained from computations with 
rational functions on the toric variety, 
many of them computer assisted. It does lead to very 
explicit formulas and parametrizations, somewhat in contrast to the 
ball quotient approach where the modular forms in question 
are hard to describe explicitly.

I'm indebted to E.\ Freitag for suggesting to undertake this study and for many discussions. I would also like to thank him and E. Carlini for assistance with the computations.

\section{Cubic surfaces their moduli space}

\subsection{} We briefly recall the basics on cubic surfaces and $E_6$, see \cite{H}
and references given there for proofs. 
We relate this to the modular orthogonal geometry used by Allcock and Freitag.

\subsection{The 27 lines.}
Any smooth cubic surface $S$ has 27 lines and there are sets of six disjoint 
lines $\{a_1,\ldots,a_6\}$. Blowing down the lines $a_i$ to points $p_i$ 
defines a birational isomorphism $S\rightarrow \PP^2$.
The images of the other 21 lines on $S$ are the 15 lines $<p_i,p_j>$ and the 6 
conics which pass through all six points except one of the $p_i$. The 
corresponding lines are denoted by $c_{ij}$ and $b_j$.
The birational inverse $\PP^2\rightarrow S$
is given by the linear system of all cubics passing through the points 
$p_1$, $\ldots$, $p_6$.

\subsection{The root system $E_6$.}
The Picard group of $S$ is isomorphic to $\ZZ^7$ and a $\ZZ$-basis is given by 
the pull-back $l$ of (the divisor class of) a line in $\PP^2$ and the classes 
of the lines $a_i$.
The intersection form is determined by 
$$
l^2=1,\quad a_i^2=-1, \quad l\cdot a_i=0,\quad a_i\cdot a_j=0
$$
for $i\neq j$. The classes of the lines are
$$
c_{ij}=l-(a_i+a_j),\qquad b_i=2l-(a_1+\ldots +\hat{a}_i+\ldots +a_6).
$$

The canonical class of $S$ is $K_S:=-3l+a_1+\ldots+a_6$ and $K_S^2=3$.
The class of a hyperplane section of $S$ is $-K_S$.
The primitive cohomology of $S$ is thus the orthogonal complement of $K_S$.
This $\ZZ$-module, 
with the bilinear form $(x,y):=-x\cdot y$, is isomorphic to the root lattice 
$Q(E_6)$ of the root system $E_6$:
$$
Q(E_6)\cong K_S^\perp :=\{x\in Pic(S):\; x\cdot K_S=0\;\}.
$$
A $\ZZ$-basis for $Q(E_6)$ is given by:
$$
\begin{array}{ccccccccc}\alpha_1=a_2-a_1,&\,& \alpha_3=a_3-a_2,&\,& 
\alpha_4=a_4-a_3,&\,&
\alpha_5=a_5-a_4,&\, &\alpha_6=a_6-a_5,\\
&&&&\phantom{x}&&&&\\
&&&&\alpha_2=l-a_4-a_5-a_6.&&&&
\end{array}
$$
This is a basis of simple roots of $E_6$:
$$
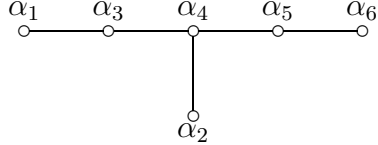
\DynkinEEE{\alpha_1}{\alpha_2}{\alpha_3}{\alpha_4}{\alpha_5}{\alpha_6}
$$
The set $E_6^+$ of positive roots of $E_6$ consists of the 36 elements in 
$Q(E_6)$ given by
$$
h_{ij}:=-a_i+a_j,\quad(i<j)\qquad h_{ijk}:=l-a_p-a_q-a_r,\quad
h:=2l-a_1-\ldots-a_6,
$$
where $\{i,j,k,l,p,q,r\}=\{1,2,\ldots,6\}$. In particular, 
$\alpha_2=h_{123}$
and with this convention our notation is compatible with that of \cite{H}.
The root system $E_6:=E_6^+\cup (-E_6^+)$ $\subset Pic(S)$ contains 
$72$ vectors, called roots.

\subsection{The Weyl group $W(E_6)$.}
The Weyl group $W(E_6)$ is the subgroup of $GL(Q(E_6))$ generated by 
reflections in the roots. We denote by $s_i$ the reflection in the hyperplane
perpendicular to the root $\alpha_i$. More generally we write $s_\alpha$,
with $\alpha\in E_6$,
for the reflection in the hyperplane perpendicular to $\alpha$.

\subsection{The orthogonal geometry.}
Allcock and Freitag use a 
non-degenerate quadratic form $Q$ on the vector space $\FF_3^5$
and its orthogonal group $O(5,3)$ to describe the combinatorics of the lines
on a cubic surfaces and of divisors on the moduli space $\cM$.
The basic facts are (\cite{Fr}, section 2):
$$
O(5,3)\cong W(E_6)\times \{\pm 1\},
$$
there are $72$ vectors with $Q(x)=-1$,
these are called the short roots (note $Q(x)=Q(-x)$). 
There are $90$ vectors with $Q(x)=-2$, the long roots, 
and there are $80$ nonzero vectors with $Q(x)=0$, 
called isotropic vectors. (See also \cite{MT}, $\S$3.)  

\subsection{Boundary divisors of ${\cM}$}\label{bdiv}
If a root is the class of an 
effective divisor on the blow up of $\PP^2$, then this effective divisor is 
a $\PP^1$ which is contracted to a node on the 
cubic surface. This sets up a correspondence between the set of irreducible 
divisors in $\cM$ parametrizing nodal cubic surfaces and $E_6^+$.
These divisors are labelled by pairs $\pm x$ of `short roots' in \cite{Fr}.

The divisor in $\cM$ corresponding to $\alpha\in E_6^+$ is denoted by $D_\alpha$
(or by $D_{ij}$ if $\alpha=h_{ij}$ etc.). These divisors are the fixed point 
sets of the corresponding reflections $s_\alpha\in W(E_6)$ in $\cM$.
The reflection $s_\alpha\in Aut(Pic(S))$ may be identified with the Picard-Lefschetz transformation associated to the general nodal cubic surface $S_0$ in $D_\alpha$.

\subsection{Lines and weights.}
Let $P(E_6)\subset Q(E_6)\otimes\QQ$ be the weight lattice of $E_6$:
$$
P(E_6):=\{x\in Q(E_6)\otimes_\ZZ\QQ:\; (x,y)\in\ZZ,\;\forall y\in Q(E_6)\,\}.
$$
The intersection number of the class $c$ of a line on $S$ with a root
is an integer, hence $c$ defines an element $x_c\in P(E_6)$. In this way one 
obtains a $W(E_6)$-orbit of $27$ weights
(which are also denoted by $a_i$, $b_i$, $c_{ij}$ with 
$1\leq i\leq 6$, $1\leq i<j\leq 6$, cf.\ \cite{H}, $\S$ 6.1.3).
Note that $a_1$ is perpendicular to 
all simple roots except $\alpha_1$ and that $(a_1,\alpha_1)=-1$, thus $a_1$ is 
minus a fundamental root of $E_6$.

\subsection{The tritangent planes and tritangent divisors.}\label{trite6}
Since hyperplane sections of $S$ correspond to cubics on the $p_i$, it is easy 
to see that there are 45 planes, the tritangent planes, which intersect $S$ in 
three lines, in Schl{\"a}fli's notation these are denoted by:
$$
(ij)=\{ a_i,b_j,c_{ij}\},\qquad
(ij.kl.mn)=\{ c_{ij},c_{kl},c_{mn}\},
$$
where$\{i,\ldots,n\}=\{1,\ldots,6\}$. Another labelling for the tritangents was given by Cayley and is used by Naruki. The dictionary between the labels is given in \cite{se2}, p.371.
The 45 tritangent divisors in $\cM$ are written as $D_t$ where $t$ is one of 
Schl\"afli's labels.
The tritangent divisors correspond to pairs $\pm x$ of long roots of \cite{Fr}.

Three lines lie in a tritangent plane iff 
the sum of their classes in $Pic(S)$ is $-K_S$ iff
the corresponding weights are linearly dependent. The orthogonal complement in $E_6$ of the span of three such weights 
is a root system of type $D_4$. If the tritangent is labelled by $t$, we will denote this $D_4$ by $t^\perp$.

\subsection{The subsystem $D_4$}\label{subd4}
An important example is the case that $t=(16)={\rm w}$. In that case $t^\perp$
is the $D_4\subset E_6$ spanned by the simple roots $\alpha_2,\;\alpha_3,\;\alpha_4$ and $\alpha_5$.
This root system is discussed in section \ref{toric}.    

\subsection{The $W(F_4)$ and tritangents.}\label{F4}
To a tritangent $t$ one associates an element 
$\gamma({\rm t})\in W(E_6)$ 
which is the product of the relections in 4 orthogonal roots in 
$t^\perp\cong D_4$. Thus $\gamma(t)$ is $-I$ on the span of $t^\perp$ and is $+I$ on the orthogonal complement which is the span on the subspace spanned by the weights corresponding to the lines in $t$.
For $t=(16)={\rm w}$ one may take 
$\gamma({\rm w})= s_2s_5s_3(s_4s_5s_3s_4)s_2(s_4s_3s_5s_4)$. 
The $\gamma(t)$'s are a conjugacy class of 45 elements in $W(E_6)$ which correspond (via their $+1$-eigenspace) with the tritangents. 
The centralizer of a $\gamma(t)$ in $W(E_6)$ is isomorphic to the Weyl group $W(F_4)$. The fixed point set of a $\gamma(t)$ on $\cC$ is the tritangent divisor $D_t$ which parametrises cubic surfaces for which the three lines in $t$ meet in one point, called an Eckart point (\cite{Naruki} \S 8).

\section{The toric variety}\label{toric}

\subsection{} For general facts on toroidal compactifications we refer to 
\cite{Fu}, for root systems see \cite{Hu}.

\subsection{The torus.}
The $D_4$-adjoint torus 
$$
T\stackrel{\cong}{\longrightarrow} (\CC^*)^4,\qquad t\longmapsto
(\lambda(t),\mu(t),\nu(t),\rho(t))
$$ 
comes with a natural identification of its character group ${\rm 
Hom}(T,\CC^*)\cong\ZZ^4$ with the sublattice
$$
M:=\langle \, e_1-e_2,\, e_2-e_3,\,e_3-e_4,\,e_3+e_4\rangle\subset 
\oplus_{i=1}^4\ZZ e_i.
$$
The lattice $M$, with the scalar product induced by the standard inner product 
on $\oplus\ZZ e_i$, is the root lattice $Q(D_4)$ of $D_4$.
We often use:
$$
{\rm Hom}(T,\CC^*)\stackrel{\cong}{\longrightarrow} M,\qquad
\lambda\mapsto e_1-e_2,\;\mu\mapsto e_3+e_4,\;\nu\mapsto e_3-e_4,\;
\rho\mapsto e_2-e_3. 
$$
For $\alpha\in M$ we define a regular function on $T$ by:
$$
f_\alpha:=\lambda^a\mu^b\nu^c\rho^d\qquad{\rm with}\quad \alpha=a(e_1-e_2)+
b(e_2-e_3)+c(e_3-e_4)+d(e_3+e_4)\in M.
$$

\subsection{The root system.}
The root system $D_4$ consists of the following 24 vectors in $M$:
$$
D_4=\{\,\pm e_i\pm e_j\;\in M:\quad 1\leq i<j\leq 4\,\}.
$$
The set 
$$
\Delta_0:=\{e_1-e_2,\, e_2-e_3,\,e_3-e_4,\,e_3+e_4\}\quad(\subset D_4)
$$ 
is a fundamental system
(or base of the root system), that is any root is a linear combination of these 
4 vector with all coefficients either positive (such a root is called positive) 
or negative. Let $N=M^*$ 
be the  dual lattice of $M$, 
$$
N:={\rm Hom}_\ZZ(M,\ZZ)=\{x\in \left(\oplus\ZZ e_i\right)^*\otimes_\ZZ\RR:\;
\langle x,\alpha \rangle\in\ZZ\quad\forall \alpha\in M\,\},
$$
here $\langle.,.\rangle$ is the pairing between $\left(\oplus\ZZ 
e_i\right)^*\otimes_\ZZ\RR$ and its dual.
Let $\{\epsilon_1,\ldots,\epsilon_4\}\subset (\oplus \ZZ e_i)^*\otimes_\ZZ\RR$ 
be the dual basis of $\{e_1,\ldots,e_4\}$. Then the basis of $N$ 
which is dual to 
$\Delta_0$ is
$$
\epsilon_1,\quad (\epsilon_1+\epsilon_2+\epsilon_3+\epsilon_4)/2,\quad
(\epsilon_1+\epsilon_2+\epsilon_3-\epsilon_4)/2, \quad \epsilon_1+\epsilon_2
\qquad (\in N).
$$

\subsection{The Weyl group.}\label{defSR}
The Weyl group $W(D_4)$ of the root system is the subgroup of $GL(M\otimes\RR)$ 
generated by the reflections in the roots (so 
$s_\alpha(\beta)=\beta-(\beta,\alpha)\alpha)$ and $(.,.)$ is the standard inner 
product on $\oplus\ZZ e_i$). This group has 192 elements and is a semidirect 
product of $S_4$ (permuting the $e_i$) and $(\ZZ/2\ZZ)^3$ (changes the sign of 
an even number of the $e_i$). 
The Weyl group acts simply transitively on the fundamental systems.

The Weyl group acts on $N$ and the 4 elements of the dual basis above are in 
distinct orbits of lengths 8, 8, 8 and 24 respectively. We define
$$
S:=\{\pm \epsilon_i\}\,\cup\,\{(\pm 
\epsilon_1\pm\epsilon_2\pm\epsilon_3\pm\epsilon_4)/2\},\qquad
R:=\{\pm\epsilon_i\pm\epsilon_j\},
$$ 
$S$ and $R$ each have 24 elements.

\subsection{The Weyl chambers.}
The (closed) Weyl chamber $C(\Delta)$ of a fundamental system 
$\Delta\;(\subset D_4)$
is the (maximal) cone in
$N\otimes_\ZZ\RR=\left(\oplus\ZZ e_i\right)^*\otimes_\ZZ\RR$ defined by:
$$
C(\Delta):=\{x\in N\otimes_\ZZ\RR\;:\quad \langle x,\alpha\rangle\geq 0\quad
\forall \alpha\in \Delta\}.
$$
If $\Delta=\{\alpha_1,\ldots,\alpha_4\}$ then the edges (i.e.\ the one 
dimensional faces) of $C(\Delta)$ are the 4 half-lines $\RR_{\geq 0}\tau_i$ 
with
$\{\tau_1,\ldots,\tau_4\}$ the dual basis of $\Delta$. The decomposition
$$
N\otimes_\ZZ\RR=\cup_{\Delta} C(\Delta)
$$
is a regular cone decomposition of the vector space $N\otimes_\ZZ\RR$,
it defines in a fan in $N$ whose faces are the faces of the 192 Weyl chambers.
This fan has 48 edges which correspond to the elements of $S\cup R$.

\subsection{The toroidal compactification.}\label{torcom}
Associated to this fan is a toric variety $\tilde{T}$,
$$
\tilde{T}=\cup_{\Delta} A(\Delta),\qquad A(\Delta)\cong \CC^4
$$
and the inclusion $T\subset A(\Delta)$ is defined by the inclusion of the rings 
of regular functions
$$
\CC[A(\Delta)]:=
\langle \;f_\alpha\;:\alpha\in M,\quad
\langle x,\alpha\rangle\geq 0\quad\forall x\in C(\Delta)\,\rangle
\;\hookrightarrow\;
\CC[T]:=\CC[\lambda^{\pm 1},\,\mu^{\pm 1},\,\nu^{\pm 1},\,\rho^{\pm 1}].
$$
For example $\CC[A(\Delta_0)]=
\CC[\lambda,\,\mu,\,\nu,\,\rho]$. 
Each edge $\RR_{\geq 0}\tau$, with $\tau\in S\cup R$, defines a divisor
$V(\tau)$ in $\tilde{T}$ (\cite{Fu}, \S 3.3) and these $48$ divisors are the 
complement of $T$ in $\tilde{T}$:
$$
\tilde{T}-T=\cup_{\tau\in S\cup R}\, V(\tau).
$$

The regular functions $f_\alpha$, $\alpha\in M$, on $T$ extend to rational 
functions on $\tilde{T}$. The divisor of $f_\alpha$ is given by:
$$
(f_\alpha)=\sum_{\tau} n_\tau V(\tau)\qquad {\rm with}\quad
n_\tau:=\langle\tau,\alpha\rangle.
$$

\subsection{Example.}\label{exal}
The divisor of $\lambda=f_{e_1-e_2}$ is given by:
$$
(\lambda)=D^+_\lambda-D^-_{\lambda}\quad{\rm with}\;
\left\{\begin{array}{lcl}
D^+_\lambda=V(\epsilon_1)+V(-\epsilon_2)+
\sum_{\pm,\pm}V((\epsilon_1-\epsilon_2
\pm\epsilon_3\pm\epsilon_4)/2)+D'\\
D^-_\lambda=V(-\epsilon_1)+V(\epsilon_2)+
\sum_{\pm,\pm}V((-\epsilon_1+\epsilon_2
\pm\epsilon_3\pm\epsilon_4)/2)+D''
\end{array} \right.
$$ 
where $D'$ and $D''$ are combinations of the divisors $V(\tau)$ with $\tau\in 
R$
with coefficients in $\{\,-2,\,-1,\,0,\,1,\,2\}$.

\subsection{The cross ratio variety.}\label{crv}
Naruki's (smooth, projective) cross ratio variety $\cC$ is obtained from the toric variety $\tilde{T}$ as follows (\cite{Naruki}, \S 10-12):
$$
\begin{array}{ccccrcccccc}
\cM&\longleftarrow&\cC&\stackrel{r}{\longleftarrow}&
\hat{T}&\stackrel{\pi''}{\longrightarrow}&
             \tT''&\stackrel{\pi'}{\longrightarrow}&
             \tT'&\stackrel{\pi_e}{\longrightarrow}&\tT.\\
\end{array}
$$

The map $\pi_e$ is the blow up of $\tT$ in the identity element $e\in T$.
The exceptional divisor $\pi_e^{-1}(e)\cong \PP^3$ is denoted by $\PP^3_{\rm w}$. The image in $\cM$ of its strict transform in $\tT''$ is the tritangent divisor $D_{\rm w}=D_{(16)}$.

The map $\pi'$ is the blow up of $\tT'$ in the strict transforms in $\tT'$
of the 12 curves in the $W(D_4)$-orbit of the curve in $\tT$ defined by $\lambda=\nu=\rho=1$. The morphism $r$ contracts the strict transforms in $\hat{T}$ of the 12 exceptional divisors in $\tT''$ to surfaces in $\cC$
and is an isomorphism on the complement (\cite{Naruki}, Prop.\ 11.3). 

The map $\pi''$ is the blow up in the strict transform 
in $\tT''$ of the 16 surfaces in the $W(D_4)$-orbit of $\mu=\rho=1$.
The 16 exceptional divisors in $\hat{T}$ map under $r$ to divisors in $\cC$, their $W(E_6)$-orbit consists of $40$ divisors, the other 24 are
the images under $r$ of the strict transforms of the $V(\tau)$'s with $\tau\in R$ (\cite{Naruki}, Prop.\ 11.2). We call these 40 divisors the cusp divisors of $\hat{T}$.

There is a morphism $\cC\rightarrow \cM$, where $\cM$ is the  moduli space of semistable marked cubic surfaces, which contracts the $40$ cusp divisors to points (cf.\ \cite{Naruki}, Introduction and \S 12), the cusps of $\cM$.  The Weyl group $W(E_6)$ acts biregularly on $\cC$ and $\cM$ and the morphism $\cC\longrightarrow \cM$ is
$W(E_6)$-equivariant.

\section{The $W(E_6)$-action on boundary divisors.}\label{e6bdiv}
 
\subsection{}
According to Naruki \cite{Naruki}, Prop.\ 11.3',
the boundary $\cC-\cM^0$ consists of two $W(E_6)$-orbits of divisors, 
one orbit is formed by the 36 boundary divisors $D_\alpha$ 
with $\alpha\in E_6^+$. The 
other orbit consists of the 40 cusp divisors and will not be of interest for us. 
In Naruki's toroidal construction, the 36 $D_\alpha$'s are parametrized
by the 12 positive roots $D_4^+$ of $D_4$ and by the 24 elements of a set of 
$S$ (see \ref{defSR}) of weights of $D_4$.
In this section we determine the corresponding $W(D_4)$-equivariant bijection 
between $E_6/\{\pm 1\}$ and $(D_4/\{\pm 1\})\cup S$, 
see table \ref{tabwe6b} for the final result.

\subsection{}\label{blT} To do the required computations, it is sufficient to 
work on the blow up of $\tT$ in the origin, rather then on  
$\cC$ or $\cM$, cf.\ \ref{crv}.   
For each positive root $\alpha\in D_4$ the closure in $\tilde{T}$
of the subtorus defined by $f_\alpha=1$ in $T$ is an irreducible divisor. Since 
it contains $e$, its pull-back to $\tilde{T}'$ has two irreducible components,
one is $\PP^3_{\rm w}$ and the other is its strict 
transform which we will denote 
by $D^1_\alpha$. The image in $\cM$ of the strict transform of $D^1_\alpha$ in $\hat{T}$ is $D_\alpha$, so these twelve divisors are labelled via $D_4={\rm w}^\perp\subset E_6$.

The other $24$ boundary divisors in $\cM$ are the images in $\cM$ 
of the strict transforms of the
$V(\beta)$ with $\beta\in S$ (\cite{Naruki}, Prop.\ 11.1). 
The Weyl group $W(D_4)$ has three orbits on $S$ and it suffices to 
identify one divisor from each orbit. That is done in the following lemma. The 
resulting labelling of all 36 divisors is given in table \ref{tabwe6b}.

\subsection{Lemma.}\label{lemact}
Let $s_1$, $s_6$ be the reflections in $W(E_6)$ defined by the roots
$\alpha_1=h_{12}$, $\alpha_6=h_{56}$ respectively. Then we have:
$$
s_1^*D^1_\lambda=V(-\epsilon_2)
$$
hence $V(-\epsilon_2)=D_{13}$. Similarly we have:
$$
\begin{array}{lcl}
s_6^*D^1_{\lambda\nu\rho}&=&V((\epsilon_1-\epsilon_2-\epsilon_3-\epsilon_4)/2),
\\
s_1^*V((\epsilon_1-\epsilon_2-\epsilon_3-\epsilon_4)/2)&=&
V((\epsilon_1+\epsilon_2+\epsilon_3+\epsilon_4)/2)
\end{array}
$$
and thus $V((\epsilon_1-\epsilon_2-\epsilon_3-\epsilon_4)/2)=D_{26}$,
$V((\epsilon_1+\epsilon_2+\epsilon_3+\epsilon_4)/2)=D_{16}$.

\ts
The divisor of the rational function $\lambda-1$ on $\tilde{T}'$ is
$$
(\lambda-1)=D^1_\lambda+\PP^3_{\rm w}-D^-_\lambda,
$$
where $D_\lambda^-$ is as in \ref{exal}. Therefore $s_1^*(\lambda-1)$ will
have exactly two effective components, one being $s_1^*D^1_\lambda$ which must 
be in the orbit of length $36$ and the other will be a tritangent divisor.
From \cite{Naruki}, p.\ 13 we have:
$$
s_1:\;\lambda\longmapsto 
\frac{\lambda\mu\nu\rho^2(1-\lambda)}{\lambda\mu\nu\rho^2-1}
$$
and hence that
$$
s_1:\;\lambda-1\longmapsto 
f_1:=\frac{1-\lambda^2\mu\nu\rho^2}{\lambda\mu\nu\rho^2-1}.
$$
Since $\lambda^2\mu\nu\rho^2=f_{2e_1}$ (note that
$\lambda^2\mu\nu\rho^2$ is not a root) and $\lambda\mu\nu\rho^2=f_{e_1+e_2}$
we see that the denominator has a pole of order one on $V(-\epsilon_2)$
but the numerator has vanishing order zero on that divisor, hence 
$V(-\epsilon_2)$ must be one of the two effective components of $(f_1)$.
The other effective component is defined by $1-\lambda^2\mu\nu\rho^2=0$,
which is the local equation of the tritangent divisor $D_{\bar{\rm x}}$ 
($\bar{\rm x}=(26)$, cf.\ Table 3 of \cite{Naruki}). 
Note that $\lambda=f_{e_1-e_2}$ and $e_1-e_2=h_{23}$, so $D^1_\lambda=D_{23}$
and that 
$s_1$ permutes the indices $1$ and $2$ of an $h_{ij}$, hence 
$s_1^*D_{23}=D_{13}$ and $s_1^*D_{(16)}=D_{(26)}$.

Using the formulas from \cite{Naruki}, p.\ 13 again we get:
$$
s_6:\;1-\lambda\nu\rho\longmapsto 
f_2:=\frac{1-\lambda\mu\nu^2\rho^2}{1-\mu\nu\rho}.
$$
Since $\lambda\mu\nu^2\rho^2=f_{e_1+e_2+e_3-e_4}$ and $\mu\nu\rho=f_{e_2+e_3}$,
we see that $V((\epsilon_1-\epsilon_2-\epsilon_3-\epsilon_4)/2)$ is 
one of the two effective components of $(f_2)$. 
The other component corresponds to the 
tritangent divisor $D_{\bar{{\rm z}}}=D_{(15)}$ defined by $1-\lambda\mu\nu^2\rho^2=0$.
Note that $\lambda\nu\rho=f_{e_1-e_4}$ and $e_1-e_4=h_{25}$, so 
$D^1_{\lambda\nu\rho}=D_{25}$
and that 
 $s_6$ permutes the indices $5$ and $6$ of an $h_{ij}$, hence 
$s_6^*D_{25}=D_{26}$ and $s_6^*D_{(16)}=D_{(15)}$.

Next we apply $s_1$ to $f_2$ and obtain:
$$
s_6:f_2\longmapsto f_3:= \frac{-\mu\rho(\lambda + \nu - \lambda\nu - 
\lambda\nu\rho - \lambda\mu\nu\rho + \lambda^2\mu\nu^2\rho^2)} 
{(\lambda\mu\rho-1)(\lambda\nu\rho-1)}
$$
In the open subset $U=A(\Delta_0)=Spec(\CC[\lambda,\mu,\nu,\rho])$,
this function is zero on $\mu=0$, which is 
$V((\epsilon_1+\epsilon_2+\epsilon_3+\epsilon_4)/2)\cap U$ 
(the zero locus in $U$ of the $i$-th element in $\{\lambda,\mu,\nu,\rho\}$ is 
the divisor corresponding to the $i$-th vector of the dual basis).
Thus we found one of the two effective components of the divisor of $f_3$.
Note that $(\rho=0)\cap U$ lies in $V(\epsilon_1+\epsilon_2)$, which is not in 
the orbit of the $36$ divisors and that the third factor of the numerator of 
$f_3$ defines the tritangent divisor labelled by $\bar{q}_1=(25)$.
Since $s_1$ permutes the indices $1$ and $2$ of an $h_{ij}$, we
get $s_1^*D_{26}=D_{16}$ and $s_1^*D_{(15)}=D_{(25)}$.
\qed

\subsection{}\label{divld}
The labelling of these 36 divisors on $\tilde{T}'$ allows us to express various 
divisors in a convenient manner. For example (cf.\ \ref{exal}):
$$
(\lambda)=D_{13}+D_{26}+D_{136}+D_{246}+D_{256}+D_{345}-
D_{12}-D_{36}-D_{126}-D_{346}-D_{356}-D_{245}
$$
and similarly:
$$
(\lambda-1)=D_{\rm w}+D_{23}-D_{12}-D_{36}-D_{126}-D_{346}-D_{356}-D_{245}.
$$

\section{The CAF-linear system.}
\subsection{}
To identify the linear system on the moduli space $\cM$ introduced by 
Coble, Allcock and Freitag and to describe the $W(E_6)$-action on it, 
we consider two divisors with support in the boundary of the 
toric variety $\tilde{T}$.

\subsection{Definition.}\label{defs}
Let $R,\;S\subset N$ be as in section \ref{defSR}.
We define divisors in $\tilde{T}$ (cf.\ \ref{torcom}) by:
$$
D_S:=\sum_{\tau\in S} V(\tau),\qquad D_R:=\sum_{\tau\in R} V(\tau).
$$

\subsection{Lemma.} We have 
$$
H^0(\tilde{T},\cO(D_S+2D_R))=\langle\,f_0=1\,\rangle\oplus \langle\, 
f_\alpha\,:
\,\alpha\in D_4\,\rangle,
$$
in particular, $\dim H^0(\tilde{T},\cO(D_S+2D_R))=25$.
The divisor $D_S+2D_R$ is very ample on $\tilde{T}$.

\ts
 The space of global sections of the line bundle
associated to a divisor $\sum n_\tau V(\tau)$ is spanned by certain 
$f_\alpha$'s:
$$
H^0(\tilde{T},\cO(\sum_\tau n_\tau V(\tau)))\,=\,\langle\,f_\alpha\,:\,
\alpha \in M\;{\rm and}\; \langle \tau , \alpha \rangle \geq -n_\tau\, \rangle 
.
$$
Thus we must find the $\alpha\in M$ with $\langle\tau,\alpha\rangle\geq -1$
for $\tau\in S$ and $\langle\tau,\alpha\rangle\geq -2$ for $\tau\in R$.
Let $\alpha=\sum m_ie_i$ with $m_i\in\ZZ$. Taking $\tau=\pm\epsilon_i\in S$ we
get $-1\leq m_i\leq 1$, taking $\tau=(\pm\epsilon_1\pm\ldots\pm\epsilon_4)/2$ 
$\in S$ we get $-2\leq \pm m_1\pm m_2\pm m_3\pm m_4\leq 2$, hence at most two 
of 
the $m_i$
are non zero and thus $\alpha=0,\,\pm e_i$ or $\pm e_i\pm e_j$ with $i\neq j$.
However $\pm e_i\not\in M$ and therefore $\alpha$ is either zero or a root.
All these $\alpha$ also satisfy $\langle\tau,\alpha\rangle\geq -2$ for $\tau\in 
R$. 

The proof of the very ampleness is standard, cf.\ \cite{Fu}, 
and since we do not really need it, we omit the proof.
\qed

\subsection{Divisors near the identity.}\label{dive}
The functions $x_r:=r-1$ with $r\in\{\lambda,\,\mu,\,\nu,\,\rho\}$
are local coordinates near the identity element $e=(1,1,1,1)\in T$.
Any rational function $f$ on $T$ which is regular in $e$ can be developed in a 
Taylor series:
$$
f=f_d+f_{d+1}+\ldots,\qquad{\rm with}\quad f_k\in 
\CC[x_\lambda,\,x_\mu,\,x_\nu,\,x_\rho]
$$
with $f_k$ homogeneous of degree $k$ and $d\geq 0$. If the polynomial
$f_d$ is not identically zero we say that $f$ vanishes to order $d$ in $e\in T$
and we write $m_e(f)=d$, $f_d$ is called the leading term of $f$.

For $\alpha=a(e_1-e_2)+\ldots+d(e_2-e_3)\in M-\{0\}$ we have: 
$$
f_\alpha-1=
(x_\lambda+1)^a(x_\mu+1)^b(x_\nu+1)^c(x_\rho+1)^d-1
=ax_\lambda+bx_\mu+cx_\nu+dx_\rho+\, H.O.T.
$$
hence $f_\alpha-1$ vanishes to order 1 at $e$ and a product $\prod_{i=1}^m 
(f_{\alpha_i}-1)$ of such functions vanishes to order $m$ at $e$.

\subsection{Definition.} \label{defV}
We define the vector space $V$ of rational functions on $\tilde{T}$ to be the 
subspace of those global sections of $\cO(D_S+2D_R)$ which vanish to order at 
least 3 at $e\in \tilde{T}$:
$$
V:=\{\,f\in\,H^0(\tilde{T},\cO(D_S+2D_R)):\; m_e(f)\geq 3\,\}.
$$

\subsection{Lemma.}\label{lemV}
The dimension of $V$ is 10. A basis for $V$, multiplied by 
$\lambda\mu\nu\rho^2$, is given in table \ref{second}.

\ts
Note that 10 is the expected dimension of $V$ since the spaces of constant, 
linear and quadratic polynomials in 4 variables have dimension $1,\,4,\,10$
respectively. Thus we only have to show that each monomial 
$x_\lambda^ax_\mu^bx_\nu^cx_\rho^d$ with $a+b+c+d\leq 2$ is the leading term of 
a function in $H^0(\tilde{T},\cO(D_S+2D_R))$. Obviously we can use $f_0=1$
to get leading term $1$ and the $r-1$ to get leading term $x_r$. 
For the roots $\alpha=\lambda\rho,\,\mu\rho,\,\nu\rho$ the leading term is
a linear combination of the leading terms of the $r-1$'s which we already have.
Subtracting these linear terms we get functions with 
the leading terms $x_\lambda x_\rho,\,x_\mu x_\rho,\,x_\nu x_\rho$.
The Taylor series of $t-1$ with $t=\lambda\nu\rho,\,\mu\nu\rho,\, 
\lambda\mu\rho$, give us, modulo the leading terms we already found,
the leading terms $x_\lambda x_\nu,\,x_\mu x_\nu,\, x_\lambda x_\nu$.
To get the $x_r^2$ use that 
$$
r^{-1}=1-x_r+x_r^2-\ldots.
$$
Thus we found all the 15 desired leading terms and we conclude that $V$ has 
codimension 15 in $H^0(\tilde{T},\cO(D_S+2D_R))$. 
\qed

\subsection{Example.}\label{f12x34x56}
 The following function lies in $V$:
$$
\nu^{-1}\rho^{-1}(\rho-1)(\lambda\nu\rho-1)(\mu\nu\rho-1)\\
=
\lambda\mu\nu\rho^2-\lambda\mu\nu\rho-\lambda\rho-\mu\rho
+\lambda+\mu
+\nu^{-1}-(\nu\rho)^{-1}.
$$
The first expression shows it vanishes to order three in $e$, 
the second that it is a linear combination of roots, 
hence it lies in $H^0(\tilde{T},\cO(D_S+2D_R))$.

\section{The action of $W(E_6)$ on the vector space $V$.}\label{we6V}

\subsection{} Naruki \cite{Naruki} defined a biregular action of $W(E_6)$ on 
$\cC$ (\cite{Naruki}, $\S5$, p.\ 13).
We show that this induces an action of $W(E_6)$ on the vector space $V$ defined 
in \ref{defV}. The vector space $V$ may be identified, via pull-back
$$
V\cong H^0(\tilde{T}',\cO(2D_R+D_S-3\PP^3_{\rm w})),
$$
where $\tilde{T}'$ is the blow up of $\tilde{T}$ in the identity element $e$ 
and 
$\PP^3_{\rm w}$ is the exceptional fiber.

The main problem is to find the images of
the divisor $D_S-3\PP^3_{\rm w}$ under $s_1,\,s_6\in W(E_6)$ and to show that
the images are linearly equivalent to this divisor. 
For this we use the following rational function:
$$ 
C_1:=\frac{(\lambda^2\mu\nu\rho^2-1)^3}
{(\lambda-1)(\lambda\rho-1)(\lambda\nu\rho-1)
(\lambda\mu\nu\rho^2-1)(\lambda\mu\nu\rho-1)(\lambda\mu\rho-1)}.
$$

\subsection{Lemma.}\label{exadiv}
The rational function $C_1$ on $\tilde{T}'$ has divisor
$$
(C_1)=3D_{\bar{\rm x}}-3\PP^3_{\rm w}+\sum_{\pm,i=2}^4 V(\pm \epsilon_i)
-D^1_{\lambda}-D^1_{\lambda\rho}-D^1_{\lambda\nu\rho}-D^1_{\lambda\mu\nu\rho^2}
-D^1_{\lambda\mu\nu\rho}-D^1_{\lambda\mu\rho}+D
$$
for some divisor $D$ which is a combination of the divisors $V(\tau)$ with 
$\tau\in R$. Here $D_{\bar{\rm x}} $ is the tritangent divisor defined by 
the strict transform of the zero locus of $\lambda^2\mu\nu\rho^2-1$ in 
$\tilde{T}$.

\ts The proof is straightforward using the formula from \ref{torcom} and the 
examples in the proof of Lemma \ref{lemact}, for example
$$
(\lambda^2\mu\nu\rho^2-1)=D_{\bar{\rm x}}+\PP^3_{\rm 
w}-2V(-\epsilon_1)-\sum_{\pm,\pm,\pm}
V((-\epsilon_1\pm\epsilon_2\pm\epsilon_3\pm\epsilon_4)/2),
$$
and the divisor of $\lambda-1$, in $\tilde{T}$, was determined in \ref{exal}.
\qed

\subsection{} For $f\in V$, the composition $f\circ s_1$ does not lie in $V$.
However, we will show that the quotient 
$(f\circ s_1)/C_1$ does lie in $V$. 
To get an action of all of $W(E_6)$ however, the
correct definition for the action of $s_1$ on $V$ is
 $s_1(f)=-(f\circ s_1)/C_1$.

\subsection{Theorem.}\label{actwe6}
The action of $W(E_6)$ on $\cC$ defines an action of $W(E_6)$ on $V$ by the 
following formulas:
$$
s_i(f):=\left\{\begin{array}{ccl}
-(f\circ s_1)/C_1&\quad& {\rm if}\;i=1,\\
f\circ s_i& &{\rm if}\;2\leq i\leq 5,\\
-(f\circ s_6)/C_6&& {\rm if}\;i=6,
\end{array}\right.
$$
here the rational maps $s_i:T\rightarrow T$ are as defined by Naruki in 
\cite{Naruki}, p.\ 13 and $C_6=C_1\circ \tau$ where 
$\tau(\lambda,\mu,\nu,\rho)=(\nu,\mu,\lambda,\rho)$.

The representation of $W(E_6)$ on $V$ is its unique 10 dimensional irreducible 
representation and is denoted by $10_s$
in \cite{Fra}.

\ts
Recall that $D_S=\sum D_\alpha$ with $\alpha\in E_6-D_4$ a positive root.
Write:
$$
D_S=D^{(0)}_S+D^{(1)}_S,\qquad D^{(0)}_S=\sum_\alpha D_\alpha
$$
where we sum over the positive roots $\alpha\in E_6$, $\alpha\not\in D_4$
which are fixed under $s_1$. Then $s_1^*D_S=D^{(0)}_S+s_1^*D^{(1)}_S$.
Since $s_1^*D_{\rm w}=D_{\bar{\rm x}}$ (cf.\ the proof of \ref{lemact}),
we get:
$$
s_1^*(D_S-3D_{\rm w})=D^{(0)}_S+s_1^*D^{(1)}_S-3D_{\bar{\rm x}}.
$$
One verifies, using the tables \ref{tabwe6} and \ref{tabwe6b} and the lemma 
above, that
$$
(C_1)=3D_{\bar{\rm x}}-3D_{\rm w}+D^{(1)}_S-s_1^*D^{(1)}_S
$$
hence $s_1^*(D_S-3D_{\rm w})+(C_1)=D_S-3D_{\rm w}$.
This suggests that 
$f\mapsto \pm(f\circ s_1)/C_1$ defines an endomorphism of $V$. 
To check this and to   
get a $W(E_6)$ representation on $V$, 
one computes matrices and checks the defining relations for $W(E_6)$ 
(we used a computer, note this direct method avoids a 
detailed discussion of the divisor $D_R$ and verifies that one has to put a `$-$' sign in the definition of $s_1$ and $s_6$). 
Since the only representations of $W(E_6)$ of dimension at most 10 are  
the trivial one, denoted by $1=1_p$, 
the 6 dimensional reflection representation $6_p$, their tensor 
products with the determinant
representation $1_n$ and $6_n$, and $10_s$, it suffices to compute the traces of a reflection 
$s_i$ (which is 0) and of a product of two commuting reflections (which has 
trace 2) to prove that $V\cong 10_s$.
\qed

\subsection{Table of a basis of $V$}\label{second}
To obtain functions in $V$, all entries have to be divided by 
$\lambda\mu\nu\rho^2$. All ten functions are in one $W(E_6)$-orbit.

{\renewcommand{\arraystretch}{2}
$$
\begin{array}{rcl}
f_1&=&(\lambda\rho-1)(\mu\rho-1)(\nu\rho-1)(\lambda\mu\nu\rho-1),\\ 
g_1&=&(\rho-1)(\lambda\mu\rho-1)(\lambda\nu\rho-1)(\mu\nu\rho-1),\\ 
f_2&=&( \mu\rho-1)( \nu\rho-1)(1 - \lambda^2\mu\nu\rho^2), \\
g_2&=&( \rho-1)( \mu\nu\rho-1)(1 - \lambda^2\mu\nu\rho^2), \\
f_3&=&( \lambda\rho-1)( \mu\rho-1)(1 - \lambda\mu\nu^2\rho^2), \\
g_3&=&( \rho-1)( \lambda\mu\rho-1)(1 - \lambda\mu\nu^2\rho^2), \\
f_4&=&\rho( \mu\rho-1)(\lambda + \nu - \lambda\nu - \lambda\nu\rho - 
\lambda\mu\nu\rho + \lambda^2\mu\nu^2\rho^2), \\
 g_4&=&\rho( \mu-1)(\lambda + \nu - \lambda\nu - \lambda\nu\rho - 
\lambda\mu\nu\rho + \lambda^2\mu\nu^2\rho^2), \\
 f_5&=&  ( \lambda\mu\rho-1)( \mu\nu\rho-1)(1 - \lambda\nu\rho^2), \\
g_5&=& ( \mu\rho-1)(\lambda\mu\nu\rho-1)(1- \lambda\nu\rho^2).
\end{array}
$$
}

\subsection{Crosses.} Allcock and Freitag construct a 10 dimensional space $W$ 
of automorphic forms on the 4-ball (\cite{Fr}, between 4.3 and 4.4)
which defines the map $\cM\hookrightarrow \PP^9$.
The vector space $W$ is spanned by certain automorphic forms which,
up to a scalar multiple, can be characterized by the fact that their divisors 
in the ball-quotient $\cM$ are crosses (\cite{Fr}, Theorem 4.6). 
A cross is defined to be a divisor 
$$
D_\alpha+D_\beta+D_\gamma+D_\delta+D_t
$$ 
where $t$ is a tritangent, defining a subroot system 
$t^\perp$ of type $D_4$ in $E_6$
(as in \ref{trite6}) and $\alpha$,$\ldots,\delta\in t^\perp\cap E_6^+$ are 
mutually perpendicular (cf.\ \cite{Fr}, Definition 3.2). 
For each tritangent $t$, there are 3 crosses
containing $D_t$, thus there are $45\cdot 3=135$ crosses. 
For example, the crosses
associated to $t=(16)$ have $\{\alpha,\ldots,\delta\}$ equal to one of the three sets:
$$
\{h_{23},\,h_{45},\,h_{123},\,h_{145}\},
\qquad
\{h_{24},\,h_{35},\,h_{124},\,h_{135}\},
\qquad
\{h_{25},\,h_{34},\,h_{125},\,h_{134}\}.
$$
The following theorem identifies $W$ with $V$ (as spaces of global sections of
a line bundle on $\cM$).

\subsection{Theorem.} \label{main}
The rational map $\tilde{F}:\tilde{T}\longrightarrow \PP^9$ defined by a basis
of the vector space $V$ defines a
$W(E_6)$-equivariant morphism
$$
F:\cC\longrightarrow \cM \;\subset \PP^9
$$
which blows down the 40 cusp divisors to the 40 cusps.
The image of $F$ is the moduli space $\cM$ which is embedded into $\PP^9$ via the map defined by Allcock and Freitag.

\ts
Using the results of $\cite{Fr}$ and the $W(E_6)$-action on $W$ and $V$,
it suffices to show that there is a function 
$f\in V\cong H^0(\tilde{T}',\cO(2D_R+D_S-3D_{\rm w}))$
such that the corresponding section has, modulo cusp divisors,
a cross as zero divisor in $\tilde{T}'$. 
In fact, the exceptional divisors in the blow ups $\pi'$ and $\pi''$ 
get blown down in the composition $\hat{T}\rightarrow \cC\rightarrow \cM$
and under push-pull via $\cM\leftarrow \hat{T}\rightarrow \tilde{T}'$
crosses in $\cM$ correspond to crosses in $\tilde{T}'$ and cusp divisors
in $\tilde{T}'$ get contracted to points in $\cM$.

Let $f=f_1$ in table \ref{second}, then the divisor in $\tilde{T}'$
of the corresponding function in $V$ is:
$$
\left({f_1\over {\lambda\mu\nu\rho^2}}\right)=
4D_{\rm w}+D_{24}+D_{124}+D_{35}+D_{135}-D_S+D'
$$
where $D'$ is a divisor with support in $D_R$.

Thus the zero divisor on $\tilde{T}'$ of the section corresponding to 
$f_1$  is $D_{\rm w}+D_{24}+D_{124}+D_{35}+D_{135}+2D_R-D'$.
Note that ${\rm w}=(16)$ and that the four roots $h_{24}$, $h_{124}$, $h_{35}$, $h_{125}$ are in $D_4=(16)^\perp$ and are perpendicular. The remaining part,
$2D_R-D'$, has support on cusp divisors. 

We observe that using the explicit bases of $V$ and the method of \cite{Fr} Corollary 7.3, one can also prove directly that $F$ factors over $\cM$ and embeds $\cM$ into $\PP^9$.
\qed

\subsection{Cross ratios.} \label{exptab}
The basis of $V$ given in \ref{second} has the property that the 
quotients $f_i/g_i$
are double ratios associated to tritangents (see the table 2 of \cite{Naruki}), and we have in fact one double ratio from each $D_4$-orbit:
$$
r({\rm w})=\frac{f_1}{g_1},\quad
r(\bar{{\rm x}})=\frac{g_2}{f_2},\quad
r(\bar{{\rm z}})=\frac{g_3}{f_3},\quad
r(\bar{{\rm q}}_1)=\frac{g_4}{f_4},\quad
r({\rm y})=\frac{f_5}{g_5}.
$$ 
(For completeness sake: ${\rm w}=(16)$, $\bar{{\rm x}}=(26)$, 
$\bar{{\rm z}}=(15)$, $\bar{{\rm q}}_1=(25)$, ${\rm y}=(16.23.45)$.) Note that the last factor in each function in \ref{second} is the local equation of the associated tritangent.

The fact that we find one cross ratio from each $D_4$ orbit 
already implies that $\cC$ is birationally isomorphic with $F(\cC)$
(use the argument of \cite{Naruki}, $\S$ 5.5).

The involution $\gamma(t)\in W(E_6)$ associated to a tritangent $t$, 
see \ref{F4},
has trace $-6$ on $V$ (cf.\ \cite{Fra}, Table II), hence it has a 2 dimensional space of invariants $V_t$ in $V$. 
There are, upto scalar multiple, 3 functions in $V_t$ whose divisors are crosses
(cf.\ \cite{Fr}, Lemma 4.5). The pairs of functions $f_i$, $g_i$ 
span such $V_t$'s. The third function in $V_{(16)}$ is:
$$
h_1:=f_1-g_1=\rho(\lambda-1)(\mu-1)(\nu-1)(\lambda\mu\nu\rho^2-1).
$$
The stabilizer $W(F_4)$ of $t$ acts on $V_t$ through the action of a 
dihedral group with $12$ elements; the subgroup $W(D_4)$ (generated by 
reflections in the long roots) acts a $S_3$ and the reflections in the short 
roots act as $-1$ on $V_t$. In fact, the elements $\sigma_1,\,\sigma_2\in 
W(F_4)$ given by Naruki in \cite{Naruki}, $\S 8$, p.\ 16 act as 
$-1$ on $V_{(16)}$.

\subsection{Complex invariants.} \label{excom}
In example \ref{f12x34x56}
we considered the following function from $V$:
$$
f=\nu^{-1}\rho^{-1}(\rho-1)(\lambda\nu\rho-1)(\mu\nu\rho-1).
$$
Its divisor satisfies, modulo components with support in $D_R$:
$$
(f)+D_S-3D_{\rm 
w}=D_{16}+D_{34}+D_{25}+D_{125}+D_{256}+D_{136}+D_{146}+D_{234}+D_{345}.
$$
The effective divisor on the right is the sum of the $D_\alpha$ where $\alpha$ 
runs over the positive roots of three mutually perpendicular $A_2$'s:
$$
\{h_{16},\;h_{125},\;h_{256}\},\qquad
\{h_{25},\;h_{234},\;h_{345}\},\qquad
\{h_{34},\;h_{136},\;h_{146}\}.
$$
There are 40 such triples of orthogonal $A_2$'s in $E_6$ which are permuted 
transitively by  $W(E_6)$ 
(\cite{H}, 6.1.5.3; this particular triple is denoted by $[16,25,34]$). 
The corresponding $40$ functions in $V$ were considered
by Coble who called them complex invariants (cf.\ \cite{Co3}, 
p.\ 340-341), see also \cite{Y}. There are $80=2\cdot 40$ functions in the $W(E_6)$-orbit of
a complex invariant, the sign of a complex invariant is not well defined.

\section{Images of divisors in $\cC$}

\subsection{} We can use Naruki's model $\cC$ and the explicit basis of $V$ to
study the moduli space $\cM\subset \PP^9$. Here we consider  various divisors in $\cM$ as subvarieties of $\PP^9$, in particular we find a nice parametrization of a tritangent divisor.

\subsection{The boundary divisors.}\label{short}

We consider the image in $\PP^9$ of one of the 36 boundary divisors 
$D_\alpha\subset \cM$ (\ref{bdiv} and section \ref{e6bdiv}). These parametrize cubic surfaces with at least one node. 
The divisor $D_\alpha$ is the fixed point set of the involution $s_\alpha$.
The trace of $s_\alpha$ on $W$ is zero, hence $W$ is the direct sum of two 
$5$-dimensional eigenspaces of $s_\alpha$. Since $F$ is equivariant for 
$W(E_6)$, $D_\alpha$ will lie in a $\PP^4$. 
The centralizer in $W(E_6)$ of the reflection $s_\alpha$ acts on the divisor $D_\alpha$ and on the eigenspaces of $s_\alpha$. This subgroup is isomorphic to $S_6$. For example if $\alpha=h$, one obtains the `standard' $S_6$ generated by all the $s_i$ except $s_2$.

In particular
we consider the image  of
$D_{345}=V(\epsilon_1)$ under $F$.
This divisor is defined by $\lambda=0$ on the open subset 
$A(\Delta_0)=Spec(\CC[\lambda,\mu,\nu,\rho])$ of $\tilde{T}$. 
Since the 10 functions listed in table \ref{second} are regular on 
$A(\Delta_0)$
and do not vanish simultaneously, we can simply take $\lambda=0$ and determine 
(the closure of) the image.
The image spans only a $\PP^4$ since the following linear functions 
vanish on this divisor (in the notation of table \ref{second}):
$$
f_1-f_2,\qquad  g_1-g_2,\qquad f_3-g_5,\qquad f_2-f_4-g_5,
\qquad g_2-g_3-f_4+g_4.
$$
The image of $\cM$ in $\PP^9$ is defined by cubics (see \cite{Fr}), and one can show that 
the image of a boundary divisor is the Segre cubic hypersurface in this $\PP^4$
(cf.\ \cite{H}, 3.2).  

\subsection{The cusp divisors.}
We consider one of the 40 cusp divisors in $\cC$ (cf.\ \ref{crv}), for example
$V(\epsilon_1+\epsilon_2)$, note that $\epsilon_1+\epsilon_2\in R$.
This divisor is defined by $\rho=0$ in $Spec(\CC[\lambda,\mu,\nu,\rho])$. 
Putting $\rho=0$ in the 10 functions in 
table \ref{second} one finds that the image of this divisor is the point
$$
(1: 1: 1: 1: 1: 1: 0: 0: 1: 1).
$$

\subsection{Tritangent divisors.}\label{trit}
The tritangent divisor $D_t$ is the fixed point set of the
involution $\gamma(t)\in W(E_6)$. Each $\gamma(t)$ has
trace $-6$ on $V$ \cite{Fra}, hence it has two eigenspaces, 
of dimension $2$ and $8$, in $V$. 
Since the dimension of the divisor $D_t$ is three 
we get $D_t\subset\PP^7$.

The centralizer of $\gamma(t)$ is isomorphic to $W(F_4)$ and this group 
acts on both $D_t$ and $\PP^7$.
We consider the case $t=(16)={\rm w}$, hence $D_t$ is birationally isomorphic
to the exceptional fiber $\PP^3_{\rm w}$ of the blow up of the torus $T$ in the identity element $e$.
Since $e$ is fixed by $W(D_4)$, we get an induced action of $W(D_4)$ on 
$\PP^3_{\rm w}$, and we will see that this action extends to a linear action of
$W(F_4)$.

The root lattice $Q(F_4)$ of $F_4$ is the lattice in $\RR^4$ generated by
the 48 roots of $F_4$ which are (cf.\ \cite{Hu}, III 12.1)
the 24 roots $\pm e_i\pm e_j$ of $D_4$ (these roots have length $2$ and are called the long roots of $F_4$) and the 24 
vectors $\pm e_i$, $(\pm e_1\pm e_2\pm e_3\pm e_4)/2$ which have length 1, the short roots of $F_4$.
$$
Q(F_4)=\langle \pm e_i\pm e_j,\;(\pm e_1\pm e_2\pm e_3\pm e_4)/2\rangle_\ZZ
\qquad(\subset\RR^4).
$$

\subsection{Theorem}\label{tritpar}
Any tritangent divisor $D_t$ is $W(F_4)$-equivariantly birationally isomorphic 
to $\PP^3$ via the map
$$
\PP^3=\PP(Q(F_4)\otimes_\ZZ\CC)\longrightarrow D_t\hookrightarrow \PP^7
$$
given by the linear system of cubics which are zero in the short roots of 
$F_4$.

\ts
Since $W(E_6)$ acts transitively on the tritangent divisors, it is sufficient
to consider the case $t=(16)$. We show that the functions from $V$ give the desired map $\PP^3_{\rm w}\rightarrow D_{(16)}$.

The local coordinate functions $\lambda-1,\ldots,\rho-1$ near $e$ induce 
projective coordinates $x_\lambda,\ldots,x_\rho$ on $\PP^3_{\rm w}$.
Since $\lambda=e_1-e_2,\ldots,\rho=e_2-e_3$ it is more convenient to use 
coordinates $y_i$ with
$$
(x_\lambda:x_\mu:x_\nu:x_\rho)=(y_1-y_2:y_3+y_4:y_3-y_4,y_2-y_3).
$$

The group $W(F_4)$ is generated by the subgroup $W(D_4)$ and 
$\sigma_1$, $\sigma_2$
given in \cite{Naruki}, $\S$8.
Using the explicit formulas for the $\sigma_i$
 one finds that these act on $\PP^3_{\rm w}$ as 
reflection in the planes $y_4=0$ and $y_1-y_2-y_3-y_4$ respectively.
(For example $\sigma_1$ interchanges $\mu$ and $\nu$ and fixes the other roots, 
thus on $\PP^3_{\rm w}$ it is the linear map which permutes 
$y_3+y_4$ and $y_3-y_4$ and fixes $y_1-y_2$ and $y_2-y_3$.)
Thus these $\sigma_i$ are reflections in the short roots.
This implies that we may identify $\PP^3_{\rm w}$ 
with $\PP(Q(F_4)\otimes_\ZZ\CC)$.

All functions in $V$ vanish to third order in $e$,
but not all vanish to fourth order,
hence restricted to $\PP^3_{\rm w}$ the map $F$ 
is given by the leading terms of third order. 
Note that $f_1$ and $g_1$ from table \ref{second} vanish to order four at $e$, 
hence the image of $\PP^3_{\rm w}$ spans at most a $\PP^7$, as we observed earlier. 
The leading terms of the other $8$ 
basis functions are cubics which all contain the 12 points:
$$
(1:0:0:0)_y,\ldots,(0:0:0:1)_y \qquad{\rm and}\quad (1:\pm 1:\pm 1:\pm 1)_y.
$$
For example,  $f_2$ from 
table \ref{second} has leading term (up to sign):
$$
(x_\mu+x_\rho)( x_\nu+x_\rho)(2x_\lambda+x_\mu+x_\nu+2x_\rho)
=(y_2+y_4)(y_2-y_4)(2y_1).
$$

One can verify that the $8$ leading cubics are independent and that
these $12$ points impose independent conditions 
on the ($20$ dimensional) space of cubics.
Thus the map $\PP^3_{\rm w}\rightarrow D_{(16)}$ is given by the subspace of cubics vanishing in these points
and the image of $\PP^3_{\rm w}$ spans a $\PP^7$.
 
We observe that the twelve basepoints are in three $D_4$-orbits of length four. Two points in 
two distinct orbits determine a line on which there is a unique line from the 
third orbit.
For example the points $(0:0:0:1)_y,\;(1:1:1:1)_y,\;(1:1:1:-1)_y$ are 
on a line. In this way we get $16$ lines, on each of 
these there are 3 base points. 
Actually, the two fourth-order leading terms 
are zero exactly on these 16 lines.
\qed

\subsection{Incidence of tritangent divisors.}\label{inctri}
The tritangent divisor $D_{(16)}$ is one of $45$ such divisors, recall that
$(16)={\rm w}=\{a_1,\,b_6,\,c_{16}\}$. The remaining $44$ 
tritangent divisors now divide into 4 groups, $44=4+4+4+2\cdot16$ as follows.
For each $l\in (16)$, there are $4$ other tritangents
containing $l$ (for example, the $\{a_1,\,b_i,\,c_{1i}\}$ for $2\leq i\leq 5$
are the other tritangents which also contain $a_1$). 

The remaining $32$ tritangents
do not have a line in common with $(16)$. These come in pairs as follows.
Given one of these $32$, say $\{l_1,l_2,l_3\}$, after a permutation of the 
indices one has that $a_1$ and $l_1$ meet 
(and $a_1$ does not meet $l_2$ and $l_3$) and thus there is a line $m_1$
such that $\{a_1,l_1,m_1\}$ is a tritangent. 
Similarly $b_6$ and $l_2$ determine
a line $m_2$ and $c_{16}$ and $l_3$ determine a line $m_3$. Now 
$\{m_1,m_2,m_3\}$ is another tritangent which has no line in common with
$(16)$. For example, $\{b_5,a_2,c_{25}\}$ determines $\{c_{15},c_{26},c_{34}\}$.
(To see all this, consider a general cubic surface and the planes 
$V_{(16)}$ and $V'$ 
spanned by the lines in $(16)$ and the $l_i$ respectively. 
These planes meet in a line which by assumption does not lie in the cubic surface. 
Thus this line meets the surface in $3$ points and through each of these points there passes exactly one line from $(16)$ and one line from the $l_i$.)

If the lines in $(16)$ all pass through an Eckart point $P_{\rm w}$ and 
similarly the lines $l_i$ all pass through an Eckart point $P'$, the line $L$
spanned by $P_{\rm w}$ and $P'$ meets the cubic surface in a third point $P''$
which is an Eckart point, being the intersection of the $m_i$.
(To see that each $m_i$ passes through $P''$, consider the plane spanned
by, say, $a_1$ and $l_1$; it cuts out $m_1$ and contains the line $L$,
hence $m_1$ meets $L$ in $P''$, similarly for the other pairs of lines.)

As a consequence, a point in the intersection of two
tritangent divisors without a common line will lie in a third 
tritangent divisor.

\subsection{Tritangent divisors and $\PP^3_{\rm w}$.}\label{leadtri}
The intersections of $\PP^3_{\rm w}$ with the other $44$ tritangent 
divisors are given by the leading terms of their equations. Those tritangent
divisors which have a line in common with ${\rm w}=\{a_1,\,b_6,\,c_{16}\}$
have a linear leading term, 
in fact one finds the following $12$ linear terms:
$$
y_i\quad(1\leq i\leq 4),\qquad{\rm and}\quad
y_1\pm y_2\pm y_3\pm y_3
$$
where the last $8$ come in two $W(D_4)$-orbits distinguished by the parity of 
the number of minus signs. For example, the tritangent $(15)=\bar{{\rm z}}$ defined by $\lambda\mu\nu^2\rho^2-1$ has leading term $y_1+y_2+y_3-y_4$. 

The tritangents which do not have a line in common with ${\rm w}$
have a leading term of 
degree two, in fact the two tritangents in a pair have the same leading term
(as they should, see the last part of \ref{inctri}).
These quadrics correspond to the 16 lines in $\PP^3_{\rm w}$ containing 
$3$ of the $12$ base points of $F$ (see \ref{trit}); each line determines a unique 
quadric by the condition that it contains the other $9$ base points (and this
quadric will not contain any of the $3$ points on the line).
For example, the tritangent divisor $\bar{q}_1=\{a_2,b_5,c_{25}\}$ is
defined by
$\lambda+\nu-\lambda\nu-\lambda\nu\rho-\lambda\mu\nu\rho+
\lambda^2\mu\nu^2\rho^2$
has leading term 
$$
Q_{25}:=y_1^2+y_2y_3-y_2y_4-y_3y_4
$$ 
which does not contain
the $3$ colinear points $(1:0:0:0)_y$, $(1:-1:-1:1)_y$ and $(1:1:1:-1)_y$.
The other tritangent divisor having the same leading term is 
$\bar{q}=\{c_{15},c_{26},c_{34}\}$ which is defined by $1-\lambda\nu\rho
-\lambda \mu \nu \rho- \lambda \mu \nu \rho^2+\lambda^2  \mu \nu \rho^2  + 
\lambda \mu \nu^2  \rho^2$.

The other $15$ quadrics can be obtained from $Q_{25}$ by the action of 
$W(D_4)$,
that is, by permuting the coordinates and changing the signs of an even number 
of the $y_i$. These quadrics are smooth and hence are isomorphic to 
$\PP^1\times\PP^1$.

\section{Equations for the moduli space.}\label{unif}

\subsection{} 
The universal marked cubic surface is embedded in a $\PP^3$-bundle over $\cC$. 
Over the moduli space $\cM^0$ of smooth marked cubic surfaces, 
this bundle is the projectivization of the tangent bundle (\cite{ACT2} $\S$ 10).
In Naruki's paper \cite{Naruki} one finds an explicit cubic 
polynomial in $R[X,Y,Z,T]$, with $R:=\CC[\lambda,\mu,\nu,\rho]$, which defines the 
universal family over an open part of $\cM$. He also gives 45 linear forms
in $R[X,Y,Z,T]$ which define the tritangent planes. 

We will verify that there are linear relations between these, 
suitably normalized, linear forms with coefficients which are elements from $V$
(note that elements from $V$ are rational functions on $T$ and thus are in the field of fractions of $R$). 
This allows us to recover the cubic 
equations found by Allcock and Freitag  which define $\cM$. 
We also find a six dimensional vector space of quintic 
polynomials, on which $W(E_6)$ acts via its standard representation, which 
define a variety $X\subset\PP^9$ whose singular locus contains $\cM$.

\subsection{}
Consider two tritangents which contain a common line. For any point in the 
interior of $\cC$, the  corresponding planes are distinct. However over one of 
the 36 boundary divisors the planes may coincide. Over a boundary 
divisor $D_\alpha$ the 6 pairs of lines  
in the double six corresponding to $\alpha$ (cf.\ \cite{H})
on the universal marked surface
specialize to the six lines through the node of the universal surface over 
$D_\alpha$.
The reflection $s_\alpha$ in $W(E_6)$ interchanges the lines in each of the six 
pairs and fixes the other 15 lines. Thus if $s_\alpha$ maps one tritangent set
to another, then the lines in the planes and thus the planes themselves will 
coincide over $D_\alpha$.

\subsection{Lemma.}
Let $t_1$, $t_2$  be two distinct tritangent sets which have a line in common. 
Then there are exactly two reflections in $W(E_6)$ which map $t_1$ to $t_2$.
The corresponding roots in $E_6$ are perpendicular.

\ts
Since $W(E_6)$ acts transitively on the set of lines, we may assume that the 
common line is $b_6$. Then the $t_i$ are of type $\{a_i,b_6,c_{i6}\}$ with 
$1\leq i\leq 5$ and applying a suitable element of $W(E_6)$ we may assume that 
$t_1=\{a_1,b_6,c_{16}\}$, $t_2=\{a_2,b_6,c_{26}\}$. 
By inspection of the lists of double sixes in \cite{H} one finds exactly one 
double six which contains the pairs $(a_1,a_2)$ and $(c_{16}, c_{26})$ (it is 
$N_{12}$) and one which contains the pairs $(a_1,c_{26})$ and $(c_{16}, a_2)$
(it is $N_{345}$).
Thus only reflections in $h_{12}$ (which permutes the indices $1$ and $2$) and 
in $h_{345}$ (which interchanges $a_1\leftrightarrow c_{26}$ and 
$a_2\leftrightarrow c_{16}$) permute these two tritangent sets.
It is easy to verify that $h_{12}$ and $h_{345}$ are perpendicular.
\qed

\subsection{}
Given three linear forms $K,\,L,\,M\in R[X,Y,Z,T]$ 
which define tritangent planes to the universal cubic surface having a 
line in common, there is a linear relation, with coefficients in $R$,
$$
AK+BL+CM=0.
$$
The next proposition shows that three tritangent planes with a line in common define three crosses. Recall that a cross is a divisor in $\cM$ 
determined by the choice of a tritangent set $t$ and on 4 perpendicular roots 
in $t^\perp\cong D_4$. 
In the example below we then verify that these crosses are the 
divisors of the coefficients in the linear relation.

\subsection{Proposition.} Let  $t_1$, $t_2$ and $t_3$ be 
tritangent sets with a line in common. Then there are crosses  $X_i$  
determined by the tritangent sets $t_i$, the pair of roots whose reflections 
interchange $t_j$ and $t_k$ (with $\{i,j,k\}=\{1,2,3\}$) and the pair of 
roots which is perpendicular to all the weights in the union of these three 
tritangent sets.

\ts
Again we use the $W(E_6)$ action, and so we may assume that 
$t_i=\{a_i,b_6,c_{i6}\}$. These span the subspace $\langle 
x_1,\,x_2,\,x_5,\,x_6\rangle$ (\cite{H}, 6.1.3) hence only the roots 
$h_{45}=-x_3+x_4$ and $h_{145}=x_3+x_4$ are perpendicular to this subspace.
The two roots whose reflections interchange $t_1$ and $t_2$
are $h_{12}$, $h_{345}$. The roots
$h_{12}$, $h_{345}$, $h_{45}$ and $h_{145}$ are orthogonal and lie in the $D_4$
perpendicular to the weights in $t_3$. Therefore there is a cross $X_3$ 
which is the sum of the tritangent divisor corresponding to $t_3$ 
and the four boundary divisors corresponding to these four roots.
Similarly one finds crosses $X_1$ and $X_2$. 
\qed

\subsection{Example.}\label{exeqq}
We consider the tritangents which contain the line $b_6$. They are:
$$
\begin{array}{cccrcr}
{\rm set}& {\rm label} &{\rm local\;equation}&&&{\rm linear\ \;form}\\
\{a_1,b_6,c_{16}\}& \;(16)={\rm w}\; &1& \qquad& &W\\
\{a_2,b_6,c_{26}\}& \;(26)=\bar{{\rm x}}\;&\lambda^2\mu\nu\rho^2 -1&
\qquad\lambda X&-&(\lambda\rho-1)(\lambda\mu\nu\rho-1)W\\
\{a_3,b_6,c_{36}\}& \;(36)={\rm x}\;&  
 \mu\nu\rho^2 -1&-X&+&(\rho-1)(\mu\nu\rho-1)W\\
\{a_4,b_6,c_{46}\}& \;(46)=x\;&-\rho(\mu\nu-1)  &X&&\\
\{a_5,b_6,c_{56}\}& \;(56)=\xi\;&\mu-\nu  &X&+&\rho(\mu-1)(\nu-1)W\\
\end{array}
$$
The conversion of the labels is given in \cite{se2},
the equation of the planes is given in \cite{Naruki}, Table 1,
but we changed the sign of $(36)$ and we multiplied the local equation
of $(46)$ by a unit.

We  write $t_i:=\{a_i,b_6,c_{16}\}$. 
Then $t_2=s_1(t_2)$, $t_3=s_3(t_2)$, $t_4=s_4(t_3)$ and $t_5=s_5(t_4)$
where $s_i$ is the reflection in $\alpha_i$.
The two roots perpendicular to the span of the
sets $t_1$, $t_2$ and $t_3$ are $h_{45}$ and $h_{123}$. 
The cross $X_1$ is then:
$$
X_1=D_{23}+D_{145}+D_{45}+D_{123}+D_{(16)}
$$
and $X_2=s_1(X_1)$, $X_3=s_3(X_2)$.

Note that $X_1$ is the divisor of the section corresponding to
$$
A_1=\rho(-1 + \lambda)(-1 + \mu)(-1 + \nu)(-1 + \lambda\mu\nu\rho^2)
(\lambda\mu\nu\rho^2)^{-1}\qquad(\in V),
$$
and  that $A_1=h_1$ in \ref{exptab}. Similarly we define
$A_2=s_1(A_1)$, $A_3=s_3(A_2)\in V$.

We define $L_{i6}\in \CC(\lambda,\ldots,\rho)[X,W]$ to be the quotient
of the linear form defining the tritangent plane $(i6)$ by the local equation of the tritangent divisor $D_{(i6)}$ as listed in the table.
One can then verify the following linear relation:
$$
A_1L_{16}+A_2L_{26}+A_3L_{36}=0.
$$

\subsection{Proposition.} \label{rank}
Let $A_i$ and $L_{ij}$ be as in Example \ref{exeqq}.
Define functions $B_i$, \ldots, $F_i\in V$ by:
$$
B_i=s_4(A_i),\quad C_i=s_3(B_i),\quad D_i=s_1(C_i),\quad E_i=s_5(D_i),
\quad F_i=s_5(B_i).
$$
Then we have $Mv=0$ where
$$
 M=\left(\begin{array}{ccccc}
A_1&A_2&A_3&0&0\\
B_1&B_2&0&B_3&0\\
C_1&0&-C_2&C_3&0\\
0&D_1&D_2& -D_3&0\\
0& E_1& E_2 & 0 & E_3\\
 F_1& F_2&0&0&- F_3\\
\end{array}\right),\qquad
v=
\left(\begin{array}{c}L_{16}\\L_{26}\\L_{36}\\L_{46}\\L_{56}\end{array}\right).
$$
In particular, $M$ has rank at most three. 

\ts
Applying the reflection $s_4$
in $\alpha_4=h_{34}$ (which permutes the indices $3$ and $4$) 
to the linear relation from Example \ref{exeqq}
we obtain a relation 
between the linear forms defining the tritangents corresponding to 
$t_1=s_4(t_1)$, $t_2=s_4(t_2)$ and $t_4=s_4(t_3)$. One verifies that this is
$B_1L_{16}+B_2L_{26}+B_3L_{46}=0$ with coefficients 
$B_i=s_4(A_i)$. Similarly, one verifies the other relations.
Since each entry of $v$ is of the form 
$a_iX+b_iW$ we see that $\ker(M)$ contains the two vectors $a=(a_1,\ldots,a_5)$ 
and $b=(b_1,\ldots,b_5)$. Thus the rank of $M$ is at 
most $5-2=3$.
 \qed

\subsection{Equations.}
To obtain equations for $\cM\subset\PP^9$ from this proposition, one chooses
a basis $X_0$,$\ldots$,$X_9$ of $V$. Then each function in $V$ is a linear form
in the $X_i$ with coefficients in $\CC$. Thus each entry of the matrix $M$
is a linear form in the $X_i$. Since the rank of $M$ is at most
$3$, the determinant of each $4\times 4$ submatrix of $M$, 
which is a degree 4 polynomial in the $X_i$, 
is identically zero as function on $\cM$.
Therefore each such determinant gives a, possibly trivial, 
quartic polynomial in the ideal of $\cM$.

\subsection{Cubics.}\label{cubrel}
To get cubic equations we consider
the following submatrix of $M$:
$$
N=\left(\begin{array}{cccc}
A_1&A_2&A_3&0\\
B_1&B_2&0&B_4\\
C_1&0&-C_3&C_4\\
\end{array}\right).
$$
The matrix $N$ has rank at most two since $Nw=0$, where $w=(L_{16},\ldots,L_{46})$,
gives two vectors in $\ker{N}$ (put $X=1$, $W=0$ and $X=0$, $W=1$ in $w$).
In particular, 
$$
\det\left(\begin{array}{ccc}
A_2&A_3&0\\
B_2&0&B_4\\
0&C_3&-C_4\\
\end{array}\right)=-A_2B_4C_3+A_3B_2C_4=0.
$$
The corresponding cubic polynomial in the $X_i$ is not identically zero in
$\CC[\ldots,X_i,\ldots]$ and is one of those found in 
 \cite{Fr} Lemma 6.3. Theorem 6.4 of that paper implies that $\cM$ is defined by the $W(E_6)$-orbit of this cubic equation.

\subsection{Quintics.}\label{quints}
One verifies that the determinant of the following submatrix of 
$M$ is a degree 5 polynomial 
in the $X_i$ which is not identically zero:
$$
M_2=\left(\begin{array}{ccccc}
A_1&A_2&A_3&0&0\\
C_1&0&-C_3&C_4&0\\
0&D_2&D_3& -D_4&0\\
0& E_2& E_3 & 0 & E_5\\
 F_1& F_2&0&0&- F_5\\
\end{array}\right).
$$
By Proposition \ref{rank} the rank of $M_2$ is at most $3$.
Therefore the determinant of any $4\times 4$ submatrix of $M_2$ is zero
on $\cM$. Since the partial deriviatives of $\det(M)$ with respect to the 
$X_i$ are linear combinations of determinants of such submatrices, we
conclude that 
the quintic hypersurface $X$ in $\PP V$ defined by $\det(M)$
is singular along moduli space of marked cubic  surfaces 
${\cal M}\subset \PP^9$.
 
Using the $10\times 5$ matrix obtained from all ${5\choose 3}=10$ linear 
relations between 3 of the 5 tritangent planes containing the line $b_6$, 
we get ${10 \choose 5}$ quintics, 
but they are either $0$ or the same as $\det(M)$ up to sign. 
It can be checked that the $W(E_6)$-orbit of such a quintic has
$27$ elements and that these quintics span a copy of 
the standard 6-dimensional representation $6_p$ of $W(E_6)$.

\section{Hunt's quintic.}\label{6dire}

\subsection{Supercrosses.}
We show how to construct 27 quintic polynomials, 
which we call supercrosses,
on $V$ which are permuted, up to sign, as the
27 lines on the cubic surface under the action of $W(E_6)$. 
We show that the supercrosses span a $6$-dimensional vector space 
on which $W(E_6)$ acts as $6_n$ 
and that they define a rational map
$$
\Sigma:\cM\longrightarrow \PP^5
$$
which maps the moduli space onto the the unique 
$W(E_6)$-invariant hypersurface of degree $5$ in $\PP^5$. 
This hypersurface was investigated by Hunt in \cite{H}.

\subsection{}
The line $a_1$ on a marked cubic surface defines a weight of $E_6$.
The roots  $\alpha_2,\ldots,\alpha_6$ are perpendicular to this weight and span
a root system, of type $D_5$, consisting of $2\cdot 20=40$ roots. 
In the notation of \cite{H}, this system is `in standard form' 
$$
a_1^\perp=\{\pm x_j\pm x_k\,:\;1\leq j<k\leq 5\}\cong D_5.
$$

Any line on a cubic surface lies in 5 tritangent planes. 
The tritangent planes containing $a_1$ are the 
$(1i)=\{a_1,b_j,c_{1j}\}$, $2\leq j\leq 6$.
The three weights corresponding to the three lines in a tritangent are linearly dependent, hence span a line, and
the orthogonal complement of the line is a root system of type $D_4$,
in fact $a_1=-(2/3)x_6$, $b_j=x_{j-1}+(1/3)x_6$, thus 
$$
\{a_1,b_j,c_{1j}\}^\perp=\{\,x_{j-1},\,x_6\}^\perp=\{\pm x_i\pm x_k:\; 
i<k,\;i,\,k\in\{1,\ldots,\widehat{j-1},\ldots,5\}\,\}\cong D_4.
$$
 
Now the main point is that the 20 positive roots which are perpendicular to 
$a_1$ split in 5 sets of 4 perpendicular roots such that each of the 5 sets
is also perpendicular to the weights corresponding to the lines in a tritangent plane containing $a_1$. Thus each line $l$ determines $5$ crosses.
In the notation of \cite{H}:
$$
\begin{array}{ccccc}
(12)=\langle a_1,\,b_2,\,c_{12}\rangle&=&\langle x_1,\,x_6\rangle
\quad&\qquad \{\pm x_2+ x_3,\;\pm 
x_4+x_5\,\}=&\quad\{h_{34},\,h_{56},\,h_{134},\,h_{156}\,\}\\
(13)=\langle a_1,\,b_3,\,c_{13}\rangle&=&\langle x_2,\,x_6\rangle
\quad&\qquad \{\pm x_1+ x_4,\;\pm 
x_3+x_5\,\}=&\quad\{h_{25},\,h_{46},\,h_{125},\,h_{146} \,\}\\
(14)=\langle a_1,\,b_4,\,c_{14}\rangle&=&\langle x_3,\,x_6\rangle
\quad&\qquad \{\pm x_1+ x_5,\;\pm 
x_2+x_4\,\}=&\quad\{h_{26},\,h_{35},\,h_{126},\,h_{135}\,\}\\
(15)=\langle a_1,\,b_5,\,c_{15}\rangle&=&\langle x_4,\,x_6\rangle
\quad&\qquad \{\pm x_1+ x_3,\;\pm 
x_2+x_5\,\}=&\quad\{h_{24},\,h_{36},\,h_{124},\,h_{136}\,\}\\
(16)=\langle a_1,\,b_6,\,c_{16}\rangle&=&\langle x_5,\,x_6\rangle
\quad&\qquad \{\pm x_1+ x_2,\;\pm 
x_3+x_4\,\}=&\quad\{h_{23},\,h_{45},\,h_{123},\,h_{145}\,\}\\
\end{array}
$$

\subsection{The functions $F_l$.}
To each cross corresponds a function, up to scalar multiple, in $V$.
Fixing one such function and applying $W(E_6)$ we find other fuctions, 
unique up to sign, whose divisors are crosses.
Fix a line $l$, then we can associate to it the function $F_l$, unique up to sign, which is the product of the 5 functions in $V$ corresponding to the $5$ crosses associated to $l$. 
The divisor of $F_l$ is then essentially the sum of the $5$ tritangent divisors
$D_t$ with $l\in t$ and the $20$ boundary divisors $D_\alpha$ with $\alpha\in l^\perp\cap E_6^+$.
If $m$ is a line and $m=\sigma(l)$ for some $\sigma\in W(E_6)$, we define
$F_m:=\det(\sigma)\sigma(F_l)$ where $\det(\sigma)$ is the determinant of $\sigma$ in the 6-dimensional reflection representation.
The $F_m$'s will be called a supercrosses, they are uniquely determined by $F_l$.

\subsection{Proposition.}
The 27 functions $F_l$ on Naruki's cross ratio variety span a 6 dimensional 
vector space. The Weylgroup $W(E_6)$ acts on this vector space as $6_n$, 
the tensor product of the standard 6 dimensional representation with its
 determinant. 

\ts
The functions $F_l$, with scalar factors suitably normalized, 
satisfy the linear relations $F_l\pm F_m\pm F_n=0$ whenever the lines $l,\,m,\,n$ are in a 
tritangent plane. From this one concludes that they span a space of dimension 6 
on which
$W(E_6)$ acts (the relations $F_{a_i}\pm F_{b_j}\pm F_{c_{ij}}=0$ 
imply one can can express the $F_{c_{ij}}$ in terms of the 
$F_{a_i}$ and $F_{b_j}$, now use the relations 
$F_{c_{ij}}\pm F_{c_{kl}}\pm F_{c_{mn}}$ to eliminate the $F_{b_j}$). 

Since reflections in the stabilizer of an $F_l$ act by as 
multiplication by $-1$ on $F_l$, 
the representation is the twist of the standard representation.
\qed

\subsection{} The theorem provides us with a $W(E_6)$-equivariant rational map
$$
\Sigma:\cM\longrightarrow \PP^5.
$$
By computing the differential of $\Sigma$ in some point of $\cM$ we found that 
it has maximal rank. Hence the (closure of the) image of $\Sigma$ is a 
$W(E_6)$-invariant hypersurface in $\PP^5$.

\subsection{Theorem.}\label{i5}
The hypersurface $\Sigma(\cM)\subset\PP^5$ is  Hunt's quintic, the unique quintic hypersurface which is $W(E_6)$-invariant. It is defined by:
$$
I_5:=\sum_l \lambda_l^5=0
$$
where $\lambda_l$ is the linear form on $\PP^5$ defined by the 
$E_6$-weight which corresponds to the line $l$.

\ts
We will show that the following sextic relation holds:
$$
\prod_{l\in A} F_l=\prod_{l\in B} F_l
$$
where  $A=\{a_1,\ldots, a_6\}$ and $B=\{b_1,\ldots,b_6\}$ form a double six
of lines.
As observed by Naruki 
(see \cite{H}, p.235), this equation is reducible, being the product of $I_5$ 
and a linear factor which is the linear form defined by the root corresponding 
to the double six given by $A$ and $B$. The $W(E_6)$-invariance of the image 
implies that the image is defined by $I_5$.

The divisors of both sides of the equation are the sum of the $6\cdot 5=30$  tritangent divisors $D_{(ij)}$ as well 
as the sum of $6\cdot 20=120$ boundary divisors.  
We already determined the 
positive roots in $a_1^\perp$ above, those in $b_1^\perp$ are:
$$
b_1^\perp=\{h_{jk}=-x_{j-1}+x_{k-1}:\;2\leq j<k\leq 6\}\cup
\{h_{pqr}:\;2\leq p<q<r\leq 6\;\}.
$$
Thus each $h_{ij}$ occurs 4 times   whereas each 
$h_{pqr}$
occurs 3 times in  both the left and the right hand side, note that $4\cdot 15+3\cdot 20=120$).
Thus, upto scalar multiple, the left and right hand side coincide.
Using the reflection $s$ in the root $h$ (note $s(a_i)=b_i$) one finds the equality.
\qed

\subsection{}
Direct computations show that the images of the 36 divisors are 36 points 
in $\PP^5$, these are the roots of $E_6$. The images of the 45 tritangent 
divisors are the 45 $\PP^3$'s in Hunt's quintic (see the proof of the theorem below).

\subsection{Theorem.}\label{thmdeg}
The rational map
$$
\Sigma:\cM\longrightarrow I_5
$$
has generic degree at least 10.

\ts  
We verified by machine computation that $\Sigma$ has maximal rank at the point
$(\lambda,\mu,\nu,\rho)=(-1,-1,2,3)\in T$. 
This point lies in the intersection of the two tritangent divisors
$(12)=\zeta$ defined by $\lambda=\mu$ and $(13)=z$ defined by $\lambda\mu=1$ (cf.\ \cite{Naruki} Table 3). These tritangents have the line $a_1$ in common.
Since $\Sigma$ is $W(E_6)$-equivariant we conclude that $\Sigma$
has maximal rank at the general point in the intersection of 
any two tritangent divisors with a line in common.

We consider the restriction of $\Sigma$ to the intersection of the tritangent divisors $D_{\rm w}=D_{(16)}$ and $D_{(26)}$ which have the line $b_6$ in common. The divisor $D_{\rm w}$ is birationally isomorphic to $\PP^3_{\rm w}$, the exceptional fiber of the blow up of $T$ in $e$, and we consider the map induced by $\Sigma$ on this $\PP^3$.
The local equation of $(26)=\bar{\rm x}$ is $\lambda^2\mu\nu\rho^2=1$ 
and its intersection with $\PP^3_{\rm w}$ is given by $y_1=0$.
(cf.\ \ref{inctri}, \ref{leadtri}). Note that $\Sigma$ has maximal rank in a general point of $\PP^3_{\rm w}\cap (y_1=0)$.

On $\PP^3_{\rm w}$ the leading terms of the $F_l$ are of degree 15 or 16 
(only for $F_{a_1}$, $F_{b_6}$ and $F_{c_{16}}$), hence the restriction of
$\Sigma$ is given by homogeneous polynomials of degree 15 and the image of $\PP^3_{\rm w}$ under $\Sigma$ lies
in the intersection of the hyperplanes defined by $a_1,\,b_6$ and $c_{16}$
which is a $\PP^3$.
After omitting leading terms which are multiples
of $y_1$ and dividing the remaining ones by their common factor $y_2y_3y_4$,
we found that $\Sigma$ restricts to $\PP^3_{\rm w}\cap (y_1=0)$ to give a map
$$
\Sigma_r:\PP^2\longrightarrow \PP^2
$$
defined by homogeneous polynomials of degree $12$. One coordinate function
is
$$
F_2:=y_3y_4(y_3-y_4)(y_3+y_4)(y_2^2-y_3y_4)^2(y_2^2+y_3y_4)^2,
$$
the other two, $F_3$ and $F_4$, are obtained by permuting the coordinates
cyclically.
All these functions satisfy
$$
F(y_2,y_3,y_4)=-(y_2y_3y_4)^8F(y_2^{-1},y_3^{-1},y_4^{-1})
$$
hence the map $\Sigma_r$ has degree at least 2.

The inverse image of a general point $(x_2:x_3:x_4)\in\PP^2$ under $\Sigma_r$
is defined by the two equations, each homogeneous of degree $12$:
$$
G_1:=x_3F_2-x_2F_3=0,\qquad G_2:=x_4F_2-x_2F_4=0.
$$
The 0-cycle defined by these equations has degree $12^2=144$, but the 
linear system defined by the $F_i$ has base points. Below we list the base 
points and their contribution to the intersection multiplicities (determined 
with computer). Here $\omega$ is a primitive cube root of unity.
$$
\begin{array}{ll}
(0:0:1),\quad (0:1:0),\quad (1:0:0),\qquad& m_P=20,\\
(0:1:\pm 1),\quad (1:0:\pm 1),\quad (1:\pm 1:0),&m_P=1,\\
(1:1:-1),\quad (1:-1:1),\quad (-1:1:1),&m_P=9,\\
(1:1:1),& m_P=9,\\
(1:\omega:\pm\omega^2),\quad (1:\omega^2:\pm \omega),\quad 
(1:-\omega:\pm\omega^2),\quad (1:-\omega^2:\pm\omega)\qquad&m_P=4.
\end{array}
$$
Thus we find that the base points contribute 
$$
3\cdot 20+6\cdot 1+3\cdot 9+1\cdot 9+8\cdot 4=134
$$
to the intersection, so there remain 10 points unaccounted for.
Since $\Sigma$ has maximal rank in a general point of this $\PP^2$,
we conclude that the degree of $\Sigma$ is at least 10.
\qed

\section{Tables.}
The following tables identify the $36$ positive roots of $E_6$, 
in the notation of Hunt \cite{H}, with the $12$ positive $D_4$ roots, 
in the notation of Naruki \cite{Naruki}, and $24$ $D_4$-weights.
We also list the functions $f_\alpha$ on $T$ corresponding to the 
positive roots $\alpha\in D_4$.

\subsection{}\label{tabwe6}
$$
\begin{array}{cclccl}
\mbox{roots of $D_4$}&f_\alpha\;&
\mbox{roots of $E_6$}\;&
\mbox{roots of $D_4$}&f_\alpha\;&
\mbox{roots of $E_6$}\\
e_1-e_2&\lambda&h_{23}=-x_1+x_2\quad&
e_1+e_2&\lambda\mu\nu\rho^2&h_{145}=x_3+x_4\\
e_1-e_3&\lambda\rho&h_{24}=-x_1+x_3&
e_1+e_3&\lambda\mu\nu\rho&h_{135}=x_2+x_4\\
e_1-e_4&\lambda\nu\rho&h_{25}=-x_1+x_4&
e_1+e_4&\lambda\mu\rho&h_{134}=x_2+x_3\\
e_2-e_3&\rho&h_{34}=-x_2+x_3&
e_2+e_3&\mu\nu\rho&h_{125}=x_1+x_4\\
e_2-e_4&\nu\rho&h_{35}=-x_2+x_4&
e_2+e_4&\mu\rho& h_{124}=x_1+x_3\\
e_3-e_4&\nu&h_{45}=-x_3+x_4&
e_3+e_4&\mu&h_{123}=x_1+x_2\\
\end{array}
$$
\subsection{}\label{tabwe6b}
$$
\begin{array}{rlrlrl}
D_4{\rm -weight}& E_6{\rm -root} & D_4{\rm -weight}& E_6{\rm -root}&D_4{-\rm 
weight}& E_6{\rm -root}\\
\epsilon_1&h_{345}&(\epsilon_1+\epsilon_2+\epsilon_3+\epsilon_4)/2 & 
h_{16}&(-\epsilon_1-\epsilon_2-\epsilon_3+\epsilon_4)/2 & h_{56} \\
\epsilon_2&h_{245} &(\epsilon_1+\epsilon_2-\epsilon_3-\epsilon_4)/2 & 
h_{236}&(-\epsilon_1-\epsilon_2+\epsilon_3-\epsilon_4)/2 & h_{46} \\
\epsilon_3&h_{235}&(\epsilon_1-\epsilon_2+\epsilon_3-\epsilon_4)/2 
&h_{246}&(-\epsilon_1+\epsilon_2-\epsilon_3-\epsilon_4)/2&h_{36}\\
\epsilon_4&h_{234}&(\epsilon_1-\epsilon_2-\epsilon_3+\epsilon_4)/2 
&h_{256}&(\epsilon_1-\epsilon_2-\epsilon_3-\epsilon_4)/2&h_{26}\\
-\epsilon_1&h_{12}&(-\epsilon_1+\epsilon_2+\epsilon_3-\epsilon_4)/2 
&h_{346}&(-\epsilon_1+\epsilon_2+\epsilon_3+\epsilon_4)/2&h_{126}\\
-\epsilon_2&h_{13}&(-\epsilon_1+\epsilon_2-\epsilon_3+\epsilon_4)/2 
&h_{356}&(\epsilon_1-\epsilon_2+\epsilon_3+\epsilon_4)/2&h_{136} \\
-\epsilon_3&h_{14}&(-\epsilon_1-\epsilon_2+\epsilon_3+\epsilon_4)/2 
&h_{456}&(\epsilon_1+\epsilon_2-\epsilon_3+\epsilon_4)/2&h_{146} \\
-\epsilon_4&h_{15}&(-\epsilon_1-\epsilon_2-\epsilon_3-\epsilon_4)/2 
&h&(\epsilon_1+\epsilon_2+\epsilon_3-\epsilon_4)/2&h_{156}\\
\end{array}
$$

\subsection{$W(E_6)$-representations.}
In the notation of Frame \cite{Fra},
the (unique) 10 dimensional representation $V$ of $W(E_6)$
is denoted by $10_s$. One has:
$$
\begin{array}{rcl}
Sym^2(10_s)&=&1+15_m+15_q+24_p,\\
Sym^3(10_s)&=&20_s+2\cdot 30_m+2\cdot 30_p+80_s,\\
Sym^4(10_s)&=&2\cdot 1+1_n+3\cdot 15_m+4 \cdot 15_q+20_p+20_s+\ldots,\\
Sym^5(10_s)&=&2\cdot 6_p+ 2\cdot 6_n+15_p+15_q+7\cdot 30_m+7\cdot 30_p+\ldots,\\
Sym^6(10_s)&=&5\cdot 1+3\cdot 1_n + 11\cdot 15_m+14\cdot 15_q+\ldots,
\end{array}
$$
here $6_p$ is the standard $6$-dimensional representation and $6_n$ is the 
tensor product of $6_p$ with its determinant. On $\PP^5$ the representations
$6_p$ and $6_n$ are the same. In particular, there are two 1-dimensional 
families
of $6$-dimensional representations in $S^5V$.

\



\begin{thebibliography}{AMR}

\bibitem[AF]{Fr} D.\ Allcock, E.\ Freitag, 
{\it Cubic surfaces and Borcherds products}, eprint math.AG/0002066.

\bibitem[ACT]{ACT2}  D.\ Allcock, J.\ A.\ Carlson, D.\ Toledo,
{\it The Complex Hyperbolic Geometry of the Moduli Space of Cubic Surfaces},
eprint math.AG/0007048.

\bibitem[C]{Co3} A.\ B.\ Coble, {\it Point sets and allied cremona groups 
(Part III)}, Trans. AMS {\bf 18}, 331--372, (1917).

\bibitem[Fr]{Fra} J.\ S.\ Frame, {\it
The classes and representations of the groups of $27$ lines and $28$ 
bitangents.} 
Ann. Mat. Pura Appl. {\bf 32}, 83--119 (1951). 


\bibitem[Fu]{Fu} W.\ Fulton, {\it Introduction to toric varieties}.
Annals of Mathematics Studies {\bf 131}, Princeton University Press (1993).


\bibitem[Hu]{Hu} J.\ E.\ Humphreys, {\it Introduction to Lie algebras and 
representation theory}. Graduate Texts
in Mathematics, Vol. {\bf 9}. Springer-Verlag, New York-Berlin, (1972).

\bibitem[H]{H} B.\ Hunt, {\it The geometry of some special arithmetic 
quotients}. Lecture Notes in
Mathematics, {\bf 1637}. Springer-Verlag, Berlin, (1996).

\bibitem[MT]{MT}  K.\ Matsumoto, T.\ Terasoma,
{\it Theta constants associated to cubic three folds}, eprint math.AG/0008024.

\bibitem[N]{Naruki} I.\ Naruki, 
{\it Cross ratio variety as a moduli space of cubic surfaces},
 Appendix by E.\ Looijenga,
Proc.\ Lond.\ Math.\ Soc. {\bf 45}, 1--30, (1982).

\bibitem[Se]{se2} J.\ Sekiguchi,
{\it The versal deformation of the $E_6$-singularity and a family of cubic 
surfaces},
J.\ Math.\ Soc.\ Japan {\bf 46}, 355--383 (1994).


\bibitem[Y]{Y} M.\ Yoshida,
{\it A W(E\_6)-equivariant projective embedding of the moduli space of cubic 
surfaces}, eprint math.AG/0002102.

\end{thebibliography}
\end{document}